\documentclass[11pt]{article} %,openbib

\usepackage[utf8x]{inputenc}
\usepackage{amsmath, amsfonts, mathrsfs, textcomp,amssymb,natbib,bbm,natbib}
\usepackage[dvips]{graphicx}

\newtheorem{theorem}{Theorem}

\newtheorem{remark}{Remark}

\newtheorem{proposition}{Proposition}
\newcommand{\ninf}{\mbox{ as } N\to\infty}

\newcommand{\bcal}{{\cal B }}
\newcommand{\ccal}{{\cal C }}
\newcommand{\dcal}{{\cal D }}

\newcommand{\norm}[1]{\left|\left|\, #1\,\right|\right|}

\newcommand{\dt}{\delta}

\newcommand{\E}{\mathrm{E}}
\newcommand{\Cov}{\mathrm{Cov}}

\newcommand{\mb}{\mathbf}

\newcommand{\beq}{ \begin{equation}}
\newcommand{\eeq}{ \end{equation}}
\newcommand{\beqr}{ \begin{eqnarray}}
\newcommand{\eeqr}{ \end{eqnarray}}
\newcommand{\beqrn}{ \begin{eqnarray*}}
\newcommand{\eeqrn}{ \end{eqnarray*}}
\newcommand{\bye}{\end{document}}

\title{Multivariate piecewise linear interpolation of a random field}
\author{Konrad Abramowicz,$\quad$ Oleg Seleznjev,\\
 Department of Mathematics and Mathematical Statistics\\
Ume{\aa } University, SE-901 87 Ume\aa , Sweden }
        \date{\today}

\begin{document}

\maketitle

\begin{abstract}
We consider a multivariate piecewise linear interpolation of a continuous random field on a $d$-dimensional cube. The approximation
performance is measured by the integrated mean square error. Multivariate piecewise linear interpolator is defined by $N$ field
observations on a locations grid (or design). We investigate the class of locally stationary random fields whose local behavior is like a
fractional Brownian field in mean square sense and find the asymptotic approximation accuracy for a sequence of designs for large $N$.
Moreover, for certain classes of continuous and continuously differentiable fields we provide the upper bound for the approximation
accuracy in the uniform mean square norm.
\end{abstract}


 \baselineskip=3.4 ex

%opening

\section{Introduction}

Let a random field $X(\mb{t}),\,\mb{t}\in[0,1]^d$, with finite second moment be observed at finite number of points.
Suppose further that the points are vertices of hyperrectangles generated by a grid in a unit hypercube.
At any unsampled point
we approximate the value of the field by a piecewise linear multivariate interpolator, which is a natural extension of a conventional one-dimensional piecewise linear
interpolator.
The approximation accuracy is measured by the integrated mean squared error.
This paper aims modelling random fields with given accuracy based on a finite number of observations.
Following \citet{Berman1974}, we extend the concept of local stationarity for random fields and focus on fields satisfying this condition.
For quadratic mean (q.m.) continuous locally stationary random fields, we derive the exact asymptotic behavior of the approximation error.
A method is proposed for determining the asymptotically optimal knot (sample points) distribution between the mesh dimensions.
We also study optimality of knot allocation along coordinates of the sampling grid.
Additionally, for q.m.\ continuous and continuously differentiable fields satisfying H\"{o}lder type conditions, we determine asymptotical upper bounds for the approximation accuracy.

The problem of random field approximation arises in many research and applied areas, like Gaussian random fields modelling  \citep{Adler2007, BrousteIstasLacroix2007}, environmental and  geosciences \citep{Christakos1992,Stein1999}, sensor networks
\citep{Zhang2005},  and image processing \citep{Pratt2007}.
The upper bound for the approximation error for isotropic random fields satisfying H\"{o}lder type conditions is given in \citet{Ritter1995}.
\citet{MullerGronbach1998} consider affine linear approximation methods and hyperbolic cross designs for fields with covariance function of tensor type.
%
%A comprehensive summary of the results for both isotropic smoothness and tensor product problems is provided in \citet{Ritter2000}.
%
An optimal allocation of the observations for Gaussian random fields with product type kernel is investigated in \citet{GronbachSchwabe1996}.
\citet{Su1997} studies limit behavior of the piecewise constant estimator for random fields with a particular form of covariance function.
\citet{Benhenni2001} investigates exact asymptotics of stationary spatial process approximation based on an equidistant sampling.
The approximation complexity and the curse of dimensionality for additive random fields are broadly discussed in \citet{Lifshits2008}.
In one-dimensional case, the piecewise linear interpolation of continuous stochastic processes is considered in, e.g., \citet{Seleznjev1996}.
Results for approximation of locally stationary processes can be found in, e.g., \citet{Seleznjev2000,HuslerPitebargSeleznjev2003,Abramowicz2011}.
\citet{Ritter2000} contains a very detailed survey of various random process and field approximation problems.
For an extensive studies of approximation problems in deterministic setting, we refer to, e.g., \citet{Nikolskii1975,deBoor2008,Kuo2009}.

The paper is organized as follows. First we introduce a basic notation.
In Section 2, we consider a piecewise multivariate linear approximation of continuous fields which local behavior is like a fractional Brownian field in mean square sense.
We derive exact asymptotics and  a formula for the optimal interdimensional knot distribution.
In the second part of this section, we provide an asymptotical upper bound for the approximation accuracy for q.m.\ continuous and differentiable fields satisfying H\"older type conditions.
In Section 3, we present the results of numerical experiments, while Section 4 contains the proofs of the statements from Section 2.

\subsection{Basic notation}

Let $X=X(\mathbf{t}),\textbf{t}\in\dcal:= [0,1]^d$, be a random field defined on a probability space $(\Omega,\mathscr{F},P)$. Assume that for every $\textbf{t}$, the random
variable $X(\textbf{t})$ lies in the normed linear space $L^2(\Omega)=L^{2}(\Omega,\mathscr{F},P)$ of
 random variables with finite second moment and identified equivalent elements with respect to $P$.
We set $||\xi||:=\left(\mathrm{E}\xi^2\right)^{1/2}$ for all $\xi \in L^2(\Omega)$ and consider the approximation
 based on the normed linear spaces of q.m.\ continuous and continuously differentiable random fields denoted by $\ccal(\dcal)$ and $\ccal^1(\dcal)$, respectively.
We define the norm for any $X\in\ccal(\dcal)$ by setting
$$
\norm{X}_p:=\left(\int_\dcal||X(\textbf{t})||^{p}d\textbf{t}\right)^{1/p},\qquad 1\leq p<\infty,
$$
and $\norm{X}_{\infty}:=\max_{\mb{t}\in\dcal}\norm{X(\mb{t})}$. For $p=2$, we call the norm \textit{integrated mean squared norm} and
the corresponding measure of approximation accuracy the \textit{integrated mean squared error} (IMSE).

Now we introduce the classes of random fields used throughout this paper.
For $k\le d$, let $\mathbf{l}=(l_1,\ldots,l_k)$ be a vector of positive integers such that $\sum_{j=1}^{k}l_j=d$, and let $L_i:=\sum_{j=1}^{i}l_j, i=0,\ldots k$, $L_0=0$, be the
sequence of its cumulative sums. Then the vector $\mathbf{l}$ defines the \mbox{\textit{l-decomposition}} of
$\dcal$ into $\dcal^1\times\dcal^2\times\ldots\dcal^k$, with the $l_j$-cube $\dcal^j=[0,1]^{l_j}$, $j=1,\ldots,k$.
For any $\mathbf{s}\in\dcal$, we denote the coordinates vector corresponding to the $j$-th component
of the decomposition by $\mathbf{s}^j$, i.e.,
$$
 \mathbf{s}^j=\mathbf{s}^j(\mb{l}):=(s_{L_{j-1}+1},\ldots, s_{L_{j}}) \in \dcal^j,\quad j=1,\ldots,k.
$$
For a vector $\boldsymbol\alpha=(\alpha_1,\ldots,\alpha_k)$, $0<\alpha_j<2$, $j=1,\ldots,k$, and the decomposition vector $\mathbf{l}=(l_1,\ldots,l_k)$,
we define
$$
\norm{\mathbf{s}}_{\boldsymbol\alpha}:=\sum_{j=1}^{k}\norm{\mathbf{s}^j}^{\alpha_j} \quad \mbox{ for all } \mathbf{s}\in\dcal
$$
%where
%  corresponds to $j$-th $l_j$ dimensional component of a vector $s$,
with the Euclidean norms $||\mathbf{s}^j||, j=1,\ldots,k$.\\

\noindent  For a random field $X\in\ccal([0,1]^d)$, we say that\\
(i) $X\in\ccal_\mathbf{l}^{\boldsymbol\alpha}([0,1]^d,C)$ if for some $\boldsymbol\alpha$,  $\mathbf{l}$,  and a positive constant $C$,
the random field $X$ satisfies the H\"{o}lder condition,  i.e.,
\beq \label{def:hcont}
    \norm{X(\mathbf{t+s})-X(\mathbf{t})}^2 \leq C \norm{\mathbf{s}}_{\boldsymbol\alpha} \qquad\mbox{ for all } \mathbf{t},\mathbf{t+s}\in[0,1]^d,
\eeq
(ii) $X\in\bcal_\mathbf{l}^{\boldsymbol\alpha}([0,1]^d,c(\cdot))$ if for some $\boldsymbol\alpha$, $\mathbf{l}$, and a vector function $c(\mb{t})=(c_1(\mb{t}),\ldots,c_k(\mb{t}))$, $\mb{t}\in[0,1]^d$, the random field $X$ is \textit{locally stationary}, i.e.,
\beq \label{def:locstat}
  \frac{\norm{X(\mathbf{t+s})-X(\mathbf{t})}^2}{\sum_{j=1}^k c_k(\mathbf{t})\norm{\mathbf{s}^j}^{\alpha_j}}\rightarrow 1\quad\mbox{as }\mathbf{s}\rightarrow 0 \mbox{ uniformly in }\mathbf{t}\in[0,1]^d,
\eeq
with positive and continuous functions $c_1(\cdot),\ldots,c_k(\cdot)$.
% \sum_{j=1}^k c_j(\mb{t})\norm{\mb{s}^j}^{\alpha_j}
We assume additionally that for $j=1,\ldots,k$, the function $c_j(\cdot)$ is invariant with respect to coordinates permutation within the $j$-th component. \\\medskip

\noindent For the classes $\ccal_\mathbf{l}^{\boldsymbol\alpha}$ and $\bcal_\mathbf{l}^{\boldsymbol\alpha}$, the withincomponent smoothness is defined by the vector $\boldsymbol\alpha=(\alpha_1,\ldots,\alpha_k)$. We denote the vector describing the smoothness for each coordinate by $\boldsymbol\alpha^*=(\alpha_1^*,\ldots,\alpha_d^*)$, where $\alpha_i^* =\alpha_j$, $i=L_{j-1}+1,\ldots,L_j$, $j=1,\ldots,k$.\medskip

\noindent\textbf{Example 1.}
Let $\mb{m}=(m_1,\ldots,m_k)$ be a decomposition vector of $[0,1]^m$, and $m=\sum_{j=1}^k m_j$.
Denote by $B_{\boldsymbol\beta,\mb{m}}(\mathbf{t})$, $\mathbf{t}\in[0,1]^m$,
$\boldsymbol\beta=(\beta_1,\ldots,\beta_k)$, $0<\beta_j<2$, $j=1,\ldots,k$,
an $m$-dimensional fractional Brownian field with covariance function
$r(\mathbf{t},\mathbf{s})=\frac{1}{2}\left(||\mathbf{t}||_{\boldsymbol\beta}+||\mathbf{s}||_{\boldsymbol\beta}-||\mathbf{t}-\mathbf{s}||_{\boldsymbol\beta}\right)$.
Then $B_{\boldsymbol\beta,\mb{m}}$ has stationary increments,
$$
||B_{\boldsymbol\beta,\mathbf{m}}(\mathbf{t+s})-B_{\boldsymbol\beta,\mathbf{m}}(\mathbf{t})||^2=||\mathbf{s}||_{\boldsymbol\beta}, \quad \mathbf{t}, \mathbf{t+s} \in[0,1]^m,
$$
and therefore, $B_{\boldsymbol\beta,\mb{m}}\in\bcal_\mathbf{m}^{\boldsymbol\beta}(\dcal,c(\cdot))$ with local stationarity functions $c_1(\mathbf{t})=\ldots=c_k(\mb{t})=1$, $\mathbf{t}\in[0,1]^m$.
In particular, if $k=1$, then $B_{\beta,m}(\mathbf{t}),\,\mb{t}\in[0,1]^m$, $0<\beta<2$, $m\in\mathbbm{N}$, is an $m$-dimensional fractal Brownian field with covariance function
\begin{equation}\label{eq:FBFdef}
r(\mathbf{t},\mathbf{s})=\frac{1}{2}\left(||\mathbf{t}||^{\beta}+||\mathbf{s}||^{\beta}-||\mathbf{t}-\mathbf{s}||^{\beta}\right),  \quad  \mathbf{t}, \mathbf{t+s} \in[0,1]^m.
\end{equation}\medskip

\noindent For $X\in\ccal^1([0,1]^d)$, we write $X_j'(\mathbf{t}),\,\mathbf{t}\in[0,1]^d$, to denote a q.m. partial derivative of $X$ with respect to the \mbox{$j$-th} coordinate, and say that
  $X\in\ccal^{1,\boldsymbol\alpha^*}([0,1]^d,{C})$ if there exist a vector $\boldsymbol\alpha^*=(\alpha_1^*,\ldots,\alpha_d^*)$
and a positive constant $C$
such that each partial derivative $X_j'$ is H\"older continuous with respect to the $j$-th coordinate, i.e., if for all $\mathbf{t},\mathbf{t+s}\in[0,1]^d$,
\begin{equation}\label{def:hdiff}
||X'_j(t_1,\ldots,t_j\!+\!s_j,\ldots,t_d)-X'_j(t_1,\ldots,t_j,\ldots,t_d)||^2\leq C |s_j|^{\alpha_j^*}\quad j=1,\ldots,d.
\end{equation}
Moreover, we say that $X\in\ccal_\mathbf{l}^{1,\boldsymbol\alpha}([0,1]^d,{C})$ with $\boldsymbol\alpha=(\alpha_1,\ldots,\alpha_k)$
 if $X\in\ccal^{1,\boldsymbol\alpha^*}([0,1]^d,{C})$ and for a given partition vector $\mathbf{l}$,  $\alpha_i:=\alpha_{L_{i-1}+1}^*=\ldots=\alpha_{L_i}^*$, $i=1,\ldots,k$.\\

Let $X$ be sampled at $N$ distinct design points $T_N$.
We consider \textit{cross regular sequences} of sampling designs $T_N:=\{\mathbf{t}_{\mathbf{i}}=(t_{1,i_1},\ldots,t_{d,i_d}): \mathbf{i}=(i_1,\ldots,i_d)$, $0\leq i_k\leq n_k^*$, $k=1,\ldots,d\}$ defined by the one-dimensional grids
$$
\int_0^{t_{j,i}}h_j^*(v)dv=\frac{i}{n_{j}^{*}},\quad i=0,1,\ldots,n_j^*,\quad j=1,\ldots,d,
$$
where $h^*_j(s)$, $s\in[0,1]$, $j=1,\ldots,d$,  are positive and continuous density functions, say, \textit{withindimensional densities}, and let
$$
{h}^*(\mb{t}):=(h_1^*(t_1),\ldots,h_d^*(t_d)).
$$
The \textit{interdimensional knot distribution} is determined by a vector function
$\pi:\mathbbm N\to\mathbbm{N}^d$:
$$
\pi^*(N):=(n_1^*(N),\ldots,n_d^*(N)),
$$
where $\lim_{N\to \infty} n_j^*(N)=\infty$, $j=1,\ldots,d$, and the condition
$$\prod_{j=1}^{d}(n_j^*(N)+1)=N$$
is satisfied. We suppress the argument $N$ for the sampling grid sizes $n_j^*=n_j^*(N)$, $j=1,\ldots,d$, when doing so causes no confusion.
Cross regular sequences are one of the possible extensions of the well known regular sequences introduced by \citet{SacksYlvisaker1966}.
The introduced classes of random fields have the same smoothness and local behavior for each coordinate of components generated by a decomposition vector $\mathbf{l}$.
Therefore in the following,  we use only approximation designs with the same within- and interdimensional knot distributions within the components.
Formally, for the partition generated by a vector $\mb{l}=(l_1,\ldots,l_k)$, we consider cross regular designs $T_N$, defined by the functions
$h:=(h_1,\ldots,h_k)$ and $\pi(N):=(n_1(N),\ldots,n_k(N))$, as follows:
$$
h_i^*(\cdot)\equiv h_j(\cdot),\quad
n_i^* = n_j,\quad i=L_{j-1}+1,\ldots,L_j,\quad j=1,\ldots,k.
$$
We call the functions $h_1(\cdot),\ldots,h_k(\cdot)$ and $\pi(N)$  \textit{withincomponent densities} and  \textit{intercomponent knot distribution}, respectively.
The corresponding property of a design $T_N$ is denoted by: $T_N$ is $cRS(h,\pi,\mb{l})$.

For a given cross regular {sampling design}, the hypercube $\dcal$ is partitioned
into $M=\prod_{j=1}^d n_j^*$ disjoint {hyperrectangles} $\dcal_{\mathbf{i}}$, $\mathbf{i}=(i_1,\ldots,i_d)$, $0\leq i_k\leq n_k^*-1$, $k=1,\ldots,d$.
Let $\mb{1}_d=(1,\ldots,1)$ denote a $d$-dimensional vector of ones. The hyperrectangle $\dcal_{\mathbf{i}}$ is
determined by the vertex
\mbox{$\mathbf{t_{i}}=(t_{1,i_1},\ldots,t_{j,i_d})$}
and the main diagonal $\mathbf{r}_{\mathbf{i}}=\mathbf{t}_{\mb{i}+\mb{1}_d}-\mathbf{t_{i}}$, i.e.,
$$
\dcal_{\mathbf{i}}:=\left\{\mathbf{t}: \mathbf{t}=\mathbf{t}_{\mathbf{i}}+\mathbf{r}_{\mathbf{i}} *\mathbf{s} , \mathbf{s}\in[0,1]^d\right\},
$$
where $'*'$ denotes the coordinatewise multiplication, i.e., for $\mathbf{x}=(x_1,\ldots,x_d)$ and $\mathbf{y}=(y_1,\ldots,y_d)$,
{$\mathbf{x}\ast\mathbf{y}:=(x_1y_1,\ldots,x_dy_d)$}.\medskip\\
For a random field $X\in{\ccal}(\dcal)$, define a
\textit{multivariate piecewise linear interpolator} (MPLI) with knots $T_N$
$$
X_N(\mathbf{t}):=X_N(X,T_N)(\mathbf{t})=\E_{\boldsymbol\eta}X(\mathbf{t_i}+\mathbf{r_i}*\boldsymbol{{\eta}}),\quad  \mathbf{t}\in\dcal_{\mathbf{i}},\,\,\mathbf{t}=\mathbf{t_{i}}+\mathbf{r_{i}}*\mathbf{s},
$$
where ${\boldsymbol \eta}=(\eta_1,\ldots,\eta_d)$ and $\eta_1,\ldots,\eta_d$ are auxiliary independent Bernoulli random variables with means $s_1,\ldots,s_d$, respectively, i.e., {$\eta_j\in {Be}(s_j)$}, $j=1,\ldots,d$.
 Such defined interpolator is continuous and piecewise linear along all coordinates.\medskip

\noindent\textbf{Example 2.} Let $d=2$, $N=4$, $\dcal=[0,1]^2$. Then $\mathbf{t}=\mathbf{s}$, $\mathbf{r}=(1,1)$,
\begin{equation*}
\begin{aligned}
X_N(\mathbf{t})=\E_{\boldsymbol \eta}X(\boldsymbol\eta)=X(0,0)(1-t_1)(1-t_2)+X(1,0)t_1(1-t_2)
+X(0,1)(1-t_1)t_2+X(1,1)t_1t_2,
\end{aligned}
\end{equation*}
and $X_N$ is a conventional bilinear interpolator \citep[see, e.g.,][]{Lancaster1986}.

We introduce some additional notation used throughout the paper. For sequences of real numbers $u_n$ and $v_n$, we write $u_n\lesssim v_n$
 if  $\lim_{n\to\infty} u_n / v_n \leq 1$ and $u_n\asymp v_n$ if there exist positive constants $c_1,c_2$ such that
$c_1 u_n\leq v_n \leq c_2 u_n$ for $n$ large enough.

\section{Results}

Let $B_{\beta,m}(\mathbf{t}),\mathbf{t}\in\mathbb{R}^{m}_{+}$, $0<\beta<2$, $m\in\mathbb{N}$, denote an $m$-dimensional fractional Brownian field with covariance function \eqref{eq:FBFdef}.
For any $\mb{u}\in\mathbb{R}^{m}_{+}$, we denote
$$
b_{\beta,m}(\mb{u}):=\int_{[0,1]^m}\norm{B_{\beta,m}(\mb{u}\ast\mb{s})-\E_{\boldsymbol\eta}B_{\beta,m}(\mb{u}\ast\boldsymbol\eta)}^2 d\mb{s},
$$
where $\boldsymbol\eta=(\eta_1,\ldots,\eta_m)$, and $\eta_1,\ldots,\eta_m$ are independent Bernoulli random variables $\eta_j\in Be(s_j),\,j=1,\ldots,m$.
Then $b_{\beta,m}(\mb{u})$ is the squared IMSE of approximation for $B_{\beta,m}(\mb{u}\ast\mb{t}),\mb{t}\in[0,1]^m$, by the MPLI
with $2^m$ observations in the vertices of unit hypercube.

In the following theorem, we provide an exact asymptotics for the IMSE of a local stationary field approximation by MPLI when a cross regular sequence of sampling designs is used.
\begin{theorem}\label{th:Main}
Let $X\in\bcal_\mathbf{l}^{\boldsymbol\alpha}(\dcal,c(\cdot))$ be a random field approximated by the MPLI $X_N(X,T_N)$, where $T_N$ is
$cRS(h,\pi,\mb{l})$. Then
$$
\norm{X-X_N}_2^2 \sim \sum_{j=1}^{k}\frac{v_j}{n_j^{\alpha_j}}>0 \mbox{ as }N \to \infty,
$$
where
\begin{equation*}
 v_j=\int_{\dcal} c_j(\mathbf{t})b_{\alpha_j,l_j}(H_j(\mb{t}^j))d\mathbf{t}>0,
\end{equation*}
and $H_j(\mb{t}^j):=(1/h_j(t_{L_{j-1}+1}),\ldots,1/h_j(t_{L_{j}}))$, $j=1,\ldots,k$.
\end{theorem}

\begin{remark}
If for the $j$-th component, the uniform withincomponent knot distribution is used, i.e., \mbox{$h_j(s)=1$}, $s\in[0,1]$, then
the asymptotic constant is reduced to
$$
v_j=\tilde b_{\alpha_j,l_j} \int_{\dcal} c_j(\mathbf{t})d\mathbf{t},
$$
where $\tilde b_{\alpha_j,l_j}:=b_{\alpha_j,l_j}(\mathbf{1}_{l_j})$.
\end{remark}

In Theorem \ref{th:Main}, the approximation accuracy is determined by the sampling grid sizes $n_j$.
The next theorem provides the asymptotically
optimal intercomponent knot distribution for a given total number of observation points $N$.
Denote by
 $$
 \rho:=\left(\sum_{i=1}^{k}\frac{l_i}{\alpha_i}\right)^{-1}\!\!\!\!=\left(\sum_{i=1}^{d}\frac{1}{\alpha_i^*}\right)^{-1}\!\!\!\!, \qquad \kappa:=\prod_{j=1}^{k}v_{j}^{l_j/\alpha_j},
 $$
where $d\!\cdot\!\rho$ is the harmonic mean of the smoothness parameters $\alpha_j^*,\,j=1,\ldots,d$.
\begin{theorem}\label{Th:DimOpt}
Let $X\in\bcal_\mathbf{l}^{\boldsymbol\alpha}(\dcal,c(\cdot))$ be a local stationary random field approximated by the MPLI $X_N(X,T_N)$, where $T_N$ is
$cRS(h,\pi,\mb{l})$. Then
\beq\label{eq:Th2}
\norm{X-X_N}_2^2 \gtrsim k\ \frac{\kappa^{\rho}}{N^{\rho}} \mbox{ as } N\to \infty.
\eeq
Moreover, for the asymptotically optimal intercomponent knot allocation,
\beq\label{opt_pi}
n_{j,opt}\sim \frac{N^{\rho/\alpha_j}\ {v_{j}^{1/\alpha_j}}}{\kappa^{\rho/\alpha_j}} \mbox{ as } N\to \infty,\quad j=1,\ldots,k,
\eeq
the equality in \eqref{eq:Th2} is attained asymptotically.
\end{theorem}
The above result agrees with the intuition that more points should be distributed in directions with lower smoothness parameters.
Note that the optimal intercomponent knot distribution leads to an increased approximation rate.
\begin{remark}
Let $X\in\bcal_\mathbf{l}^{\boldsymbol\alpha}(\dcal,c(\cdot))$ with $k=d$ and $\alpha_i\neq\alpha_j$ for some $i,j=1,\ldots,d$, and  $\underline{\alpha}:=\min_{i=1,\ldots,d}\alpha_i$, i.e.,
$\rho>\underline{\alpha}$.
 Consider the approximation with uniform
intercomponent knot distribution, $n_1=\cdots=n_d\sim N^{1/d}$. Then by Theorem \ref{th:Main}, we have
$$
\norm{X-X_N}_2\asymp \frac{1}{N^{\underline{\alpha}/{(2d)}}}.
$$
On the other hand, the sampling distribution \eqref{opt_pi} gives
$$
\norm{X-X_N}_2\asymp \frac{1}{N^{\rho/2}}<\frac{1}{N^{\underline{\alpha}/{(2d)}}}.
$$
\end{remark}

\noindent\textbf{Example 3.} Let $d=k=2$, $\alpha_1=2/3$, $\alpha_2=5/3$. Then for $n_1=n_2$,  the approximation rate is $N^{-\underline\alpha/2d}=N^{-1/6}$ while using the asymptotically optimal intercomponent distribution we obtain
the rate $N^{-\rho/2}=N^{-1/4.2}<N^{-1/6}$.\bigskip

In general setting, numerical
procedures can be used for finding optimal densities. However, in practice such methods are very computationally demanding.
We present a simplification of the asymptotic constant expression for one-dimensional components.
Further, in this case, we provide the exact formula for the density minimizing the asymptotic constant. For a random field $X\in\bcal_\mb{l}^{\boldsymbol\alpha}(\dcal,c(\cdot))$, define the \textit{integrated local stationarity functions}
$$
\begin{aligned}
C_j(t_{L_j})&:=\int_{[0,1]^{d-1}}c_j(\mathbf{t}) dt_1 \ldots dt_{{L_j-1}}dt_{{L_j}+1}\ldots dt_{d},\quad &t_{L_j}\in[0,1],\quad j=1,\ldots,k.\\
\end{aligned}
$$
Moreover, for $0<\beta<2$, let
$$
a_{\beta}:=\frac{2}{(\beta+1)(\beta+2)}-\frac{1}{6}.
$$

\begin{proposition}\label{Prop:OneDim} Let $X\in\bcal_\mb{l}^{\boldsymbol\alpha}(\dcal,c(\cdot))$ be a random field approximated by the MPLI $X_N(X,T_N)$, where $T_N$ is $cRS(h,\pi,\mb{l})$.
 If for some $j$, $1\leq j\leq k$,  $l_j=1$, then for any regular density $h_j(\cdot)$, we have
 \begin{equation*}
 v_j=a_{\alpha_j}\int_0^1 C_j(t_{L_j}) h_j(t_{L_j})^{-\alpha_j} dt_{L_j}.
 \end{equation*}
The density minimizing $v_j$ is given by
\begin{equation*}
h_{j,opt}(t_{L_j})=\frac{C_j(t_{L_j})^{\gamma_j}}{\int_0^1 C_j(\tau_{L_j})^{\gamma_j}d\tau_{L_j}}, \qquad t_{L_j}\in[0,1],
\end{equation*}
where $\gamma_j:=1/(1+\alpha_j)$. Furthermore, for such chosen density, we get $$v_{j,opt}=a_{\alpha_j}\norm{C_j}_{\gamma_j}.$$
\end{proposition}

In the subsequent proposition, we give an upper bound for the approximation error
together with expressions for generating densities minimizing this upper bound, called \textit{suboptimal} densities.

\begin{proposition}\label{th:Suboptimal}
Let $X\in\bcal_\mathbf{l}^{\boldsymbol\alpha}(\dcal,c(\cdot))$ be a random field approximated by the MPLI $X_N(X,T_N)$, where $T_N$ is
$cRS(h,\pi,l)$. Then
$$
\norm{X-X_N}_2^2\lesssim \sum_{j=1}^{k}\frac{w_j}{n_j^{\alpha_j}} \ninf,
$$
where
\begin{equation*}
w_j= {l_j^{1+\alpha_j/2}}\left(a_{\alpha_j}+\frac{1}{6}\right) \int_0^1 C_j(t_{L_j}) h_j(t_{L_j})^{-\alpha_j}dt,\quad j=1,\ldots,k.
\end{equation*}
The density minimizing $w_j$ is given by
\begin{equation*}
h_{j,subopt}(t_{L_j})=\frac{C_j(t_{L_j})^{\gamma_j}}{\int_0^1 C_j(\tau_{L_j})^{\gamma_j}d\tau_{L_j}}, \qquad t_{L_j}\in[0,1],
\end{equation*}
where $\gamma_j:=1/(1+\alpha_j)$, $j=1,\ldots,k$. Furthermore, for such chosen densities, we get
$$w_{j,subopt}={l_j^{1+\alpha_j/2}}\left(a_{\alpha_j}+\frac{1}{6}\right)\norm{C_j}_{\gamma_j},\quad j=1,\ldots,k.$$
\end{proposition}

Now we focus on random fields satisfying the introduced H\"{o}lder type conditions. In this case, we provide results for the \textit{uniform mean square norm} of
approximation error $\norm{X-X_N}_\infty$.
The following proposition provides an upper bound for the accuracy of MPLI for H\"{o}lder classes of continuous and continuously
differentiable fields.
\begin{proposition}\label{th:HoldType}
Let $X\in\ccal(\dcal)$ be a random field approximated by the MPLI $X_N(X,T_N)$, where $T_N$ is
$cRS(h,\pi,\mb{l})$.
\begin{itemize}
 \item[(i)] If $X\in\ccal_\mb{l}^{\boldsymbol\alpha}(\dcal,C)$,
       then
\begin{equation}\label{eq:HoldType1}
\norm{X-X_N}_\infty\leq \sqrt{C}\sum_{j=1}^{k}\frac{c_j}{n_j^{\alpha_j/2}}
\end{equation}
for positive constants $c_1,\ldots,c_k$.
 \item[(ii)] If $X\in\ccal_\mb{l}^{1,\boldsymbol\alpha}(\dcal,C)$,
then
\begin{equation}\label{eq:HoldType2}
\norm{X-X_N}_\infty\leq \sqrt{C} \sum_{j=1}^{k}\frac{d_j}{n_j^{1+\alpha_j/2}}
\end{equation}
for positive constants $d_1,\ldots,d_k$.
\end{itemize}
\end{proposition}
\begin{remark} It follows from the proof of Proposition \ref{th:HoldType} that \eqref{eq:HoldType1} holds if
$$
c_j^2=2^{-\alpha_j} l_j^{1+\alpha_j/2} D_j^{\alpha_j},\quad j=1,\ldots,k,
$$
where $D_j:=1/\min_{s\in[0,1]}h_j(s)$, $j=1,\ldots,k$. Therefore the constants depend only on the parameters of the H\"{o}lder class and the corresponding sampling design.
Similar formulas can be obtained for $d_1,\ldots,d_k$ in \eqref{eq:HoldType2}.
\end{remark}
In addition, we provide the intercomponent knot distribution leading to an increased rate of the upper bounds obtained in Proposition \ref{th:HoldType}.
\begin{remark}\label{rem:HoldOpt} Let $X\in\ccal(\dcal)$ be a random field approximated by the MPLI $X_N(X,T_N)$, where $T_N$ is
$cRS(h,\pi,l)$.
\begin{itemize}
 \item[(i)] If $X\in\ccal_\mb{l}^{\boldsymbol\alpha}(\dcal,C)$ and
$n_j\sim N^{\rho_0/\alpha_j},j=1,\ldots,k,$
where $\rho_0=(\sum_{i=1}^{k}l_i/\alpha_i)^{-1}$, then
$$
\norm{X-X_N}_\infty=O(N^{-\rho_0/2}) \ninf.
$$
 \item[(ii)] If $X\in\ccal_{\mb{l}}^{1,\boldsymbol\alpha}(\dcal,C)$ and
$n_j\sim N^{\rho_1/(2+\alpha_j)},j=1,\ldots,k,$
where $\rho_1=(\sum_{i=1}^{k}l_i/(2+\alpha_i))^{-1}$, then
$$
\norm{X-X_N}_\infty=O(N^{-\rho_1/2}) \ninf.
$$
\end{itemize}
\end{remark}
The approximation rates obtained in the above remark are optimal in a certain sense, i.e., the rate of
convergence can not be improved in general for random fields satisfying H\"{o}lder type condition \citep[see, e.g.,][]{Ritter2000}.
Moreover, these rates correspond to the optimal approximation rates for anisotropic Nikolskii-H\"{o}lder classes \citep[see, e.g., ][]{Yanjie2000}, which are
deterministic analogues of the introduced H\"{o}lder classes.

\section{Numerical Experiments}
In this section, we present some examples illustrating the obtained results.
For given  knot densities and covariance functions,  first the pointwise approximation errors are found  analytically.
Then numerical integration is used to evaluate the approximation errors on the entire unit hypercube.
Let
$$
\dt_N(h,\pi)(\mb{t})=\dt_N(X,X_N,T_N(h,\pi,\mathbf{l}))(\mathbf{t}):=X(\mathbf{t})-X_N(X,T_N(h,\pi,\mathbf{l}))(\mathbf{t}),\qquad \mathbf{t}\in[0,1]^d,
$$
be the deviation field for the approximation of $X$ by the MPLI with $N$ knots, where $T_N$ is $cRS(h,\pi,\mb{l})$,
and write
$$
e_N(h,\pi):=\norm{\dt_N(h,\pi)}_2
$$
for the corresponding IMSE.
We write $h_{uni}(\cdot)$, to denote the vector of withincomponent uniform densities. Analogously, by $\pi_{uni}(\cdot)$ we denote the uniform interdimensional knot distribution, i.e.,
$n_1=\ldots=n_k$.

\begin{figure}[hbt]
\begin{center}
\includegraphics[height=1.8in]{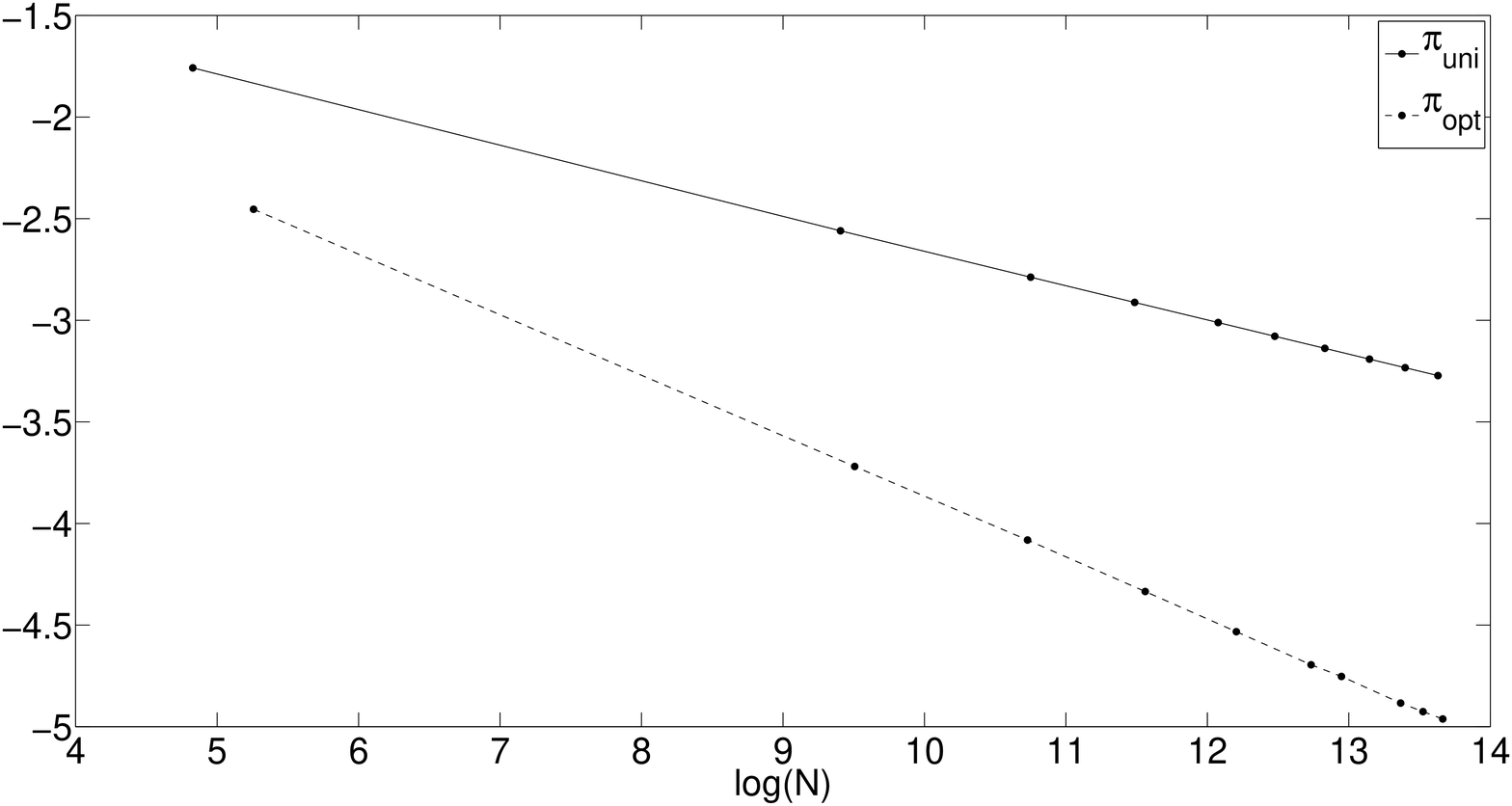}
\end{center}
\caption{The (fitted) plots of $e_N^2(h_{uni},\pi_{uni})$ (solid line), $e_N^2(h_{uni},\pi_{opt})$ (dash line) versus $N$ in a log-log scale.}
\label{fg:TwoFBM}
\end{figure}

\noindent\textbf{Example 4.}
Let $\dcal=[0,1]^3$ and
$$
X(\mb{t})=B_{\boldsymbol\alpha,\mb{l}}(\mathbf{t}),
$$
where $\boldsymbol\alpha=(1/2,3/2)$ and $\mathbf{l}=(1,2)$. Then
$X\in\bcal_{\mb{l}}^{\boldsymbol\alpha}([0,1]^3,c(\cdot))$, with $c(\mb{t})=(1,1),\,\mb{t}\in[0,1]^3$, $k=2$, $\boldsymbol\alpha^*=(1/2,3/2,3/2)$.
We compare behavior of $e_N(h_{uni},\pi_{uni})$ and $e_N(h_{uni},\pi_{opt})$, where $\pi_{opt}$ given by Theorem~\ref{Th:DimOpt}.
Observe that by using the asymptotically optimal intercomponent distribution, we obtain a gain in the rate of approximation.
Figure \ref{fg:TwoFBM} shows
the (fitted) values of the squared IMSEs $e_N^2(h_{uni},\pi_{uni})$ and $e_N^2(h_{uni},\pi_{opt})$ in a log-log scale. In such scale, the slopes of fitted lines correspond to the rates of approximation.
These plots represent the following asymptotic behavior:
$$
\begin{aligned}
e_N^2(h_{uni},\pi_{uni}) &\sim  0.3667 N^{-1/6} +0.0935 N^{-1/2}\sim 0.3667 N^{-1/6},&\\
e_N^2(h_{uni},\pi_{opt}) &\sim  0.4245 N^{-3/10} &\ninf.
\end{aligned}
$$

\begin{figure}[hbt]
\begin{center}
\begin{tabular}{cc}
\includegraphics[height=1.8in]{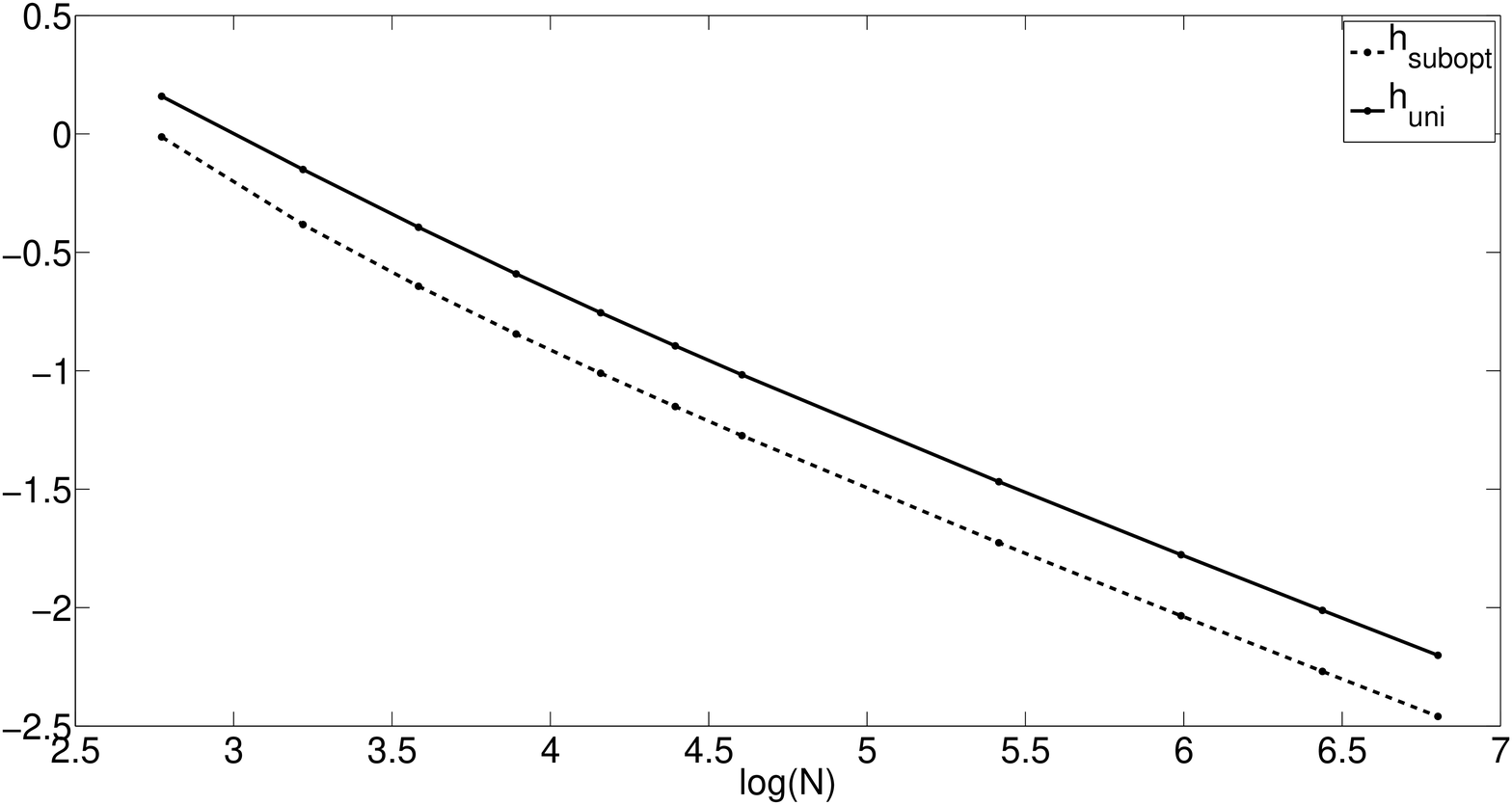} & \includegraphics[height=1.8in]{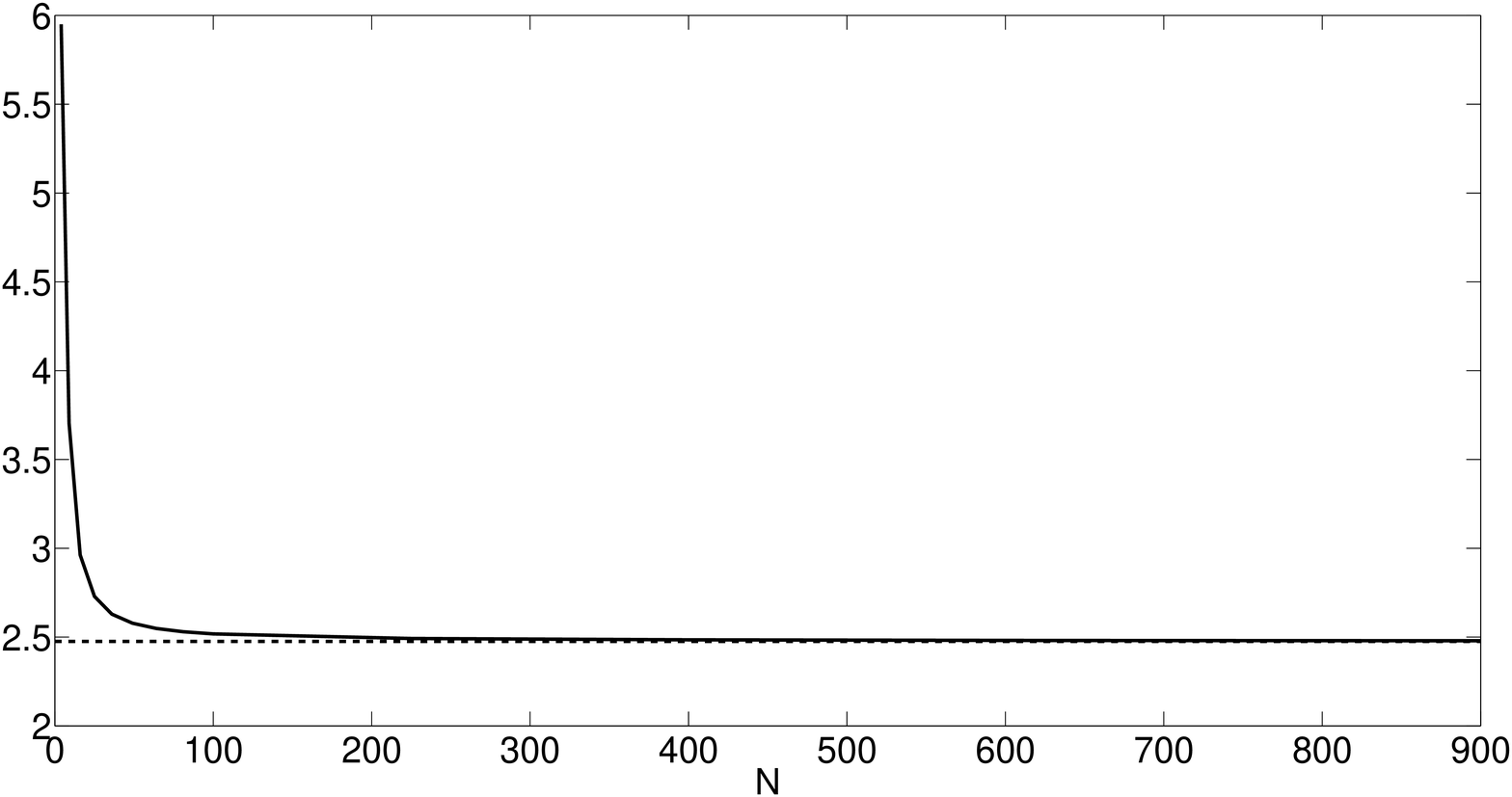}\\
\small{(a)}&
\small{(b)}
\end{tabular}
\end{center}
\caption{(a) The (fitted) plots of $e_N^2(h_{uni},\pi_{uni})$ (dashed line) and $e_N^2(h_{subopt},\pi_{uni})$ (solid line) versus $N$ in a log-log scale.
(b) The convergence of $N^{0.5}e_N^2(h_{subopt},\pi_{uni})$ (solid line) to the asymptotic constant (dashed line).}
\label{fg:TwoDimIso}
\end{figure}

\noindent\textbf{Example 5.}
Let $\dcal=[0,1]^2$ and define $X(\mathbf{t})=X(t_1,t_2)$ to be a zero mean Gaussian field with covariance function
$$
\Cov(X(\mathbf{t}),X(\mathbf{s}))=\frac{1}{(||\mathbf{t}||^2+0.1)}\frac{1}{(||\mathbf{s}||^2+0.1)}\exp(-||\mathbf{t-s}||).
$$
Then $X\in\bcal_{\mathbf{l}}^{\boldsymbol\alpha}([0,1]^2,c(\cdot))$ with $c(\mathbf{t})=c_1(\mb{t})=2/(||\mathbf{t}||^2+0.1)^2 ,\,\mathbf{t}\in[0,1]^2$, $\boldsymbol\alpha=1$,
$\boldsymbol\alpha^*=(1,1)$, $\mathbf{l}=2$, and $k=1$.
The field has one component, hence the uniform interdimensional knot distribution is used.
Theorem \ref{th:Suboptimal} provides the formula for the suboptimal withincomponent density.
Figure \ref{fg:TwoDimIso}(a) shows the (fitted) values of the squared IMSEs $e_N^2(h_{uni},\pi_{uni})$ and $e_N^2(h_{subopt},\pi_{uni})$.
Figure \ref{fg:TwoDimIso}(b) demonstrates the convergence of the scaled squared approximation error $N^{0.5}e_N^2(h_{subopt},\pi_{uni})$ to the asymptotic
constant obtained in Theorem~\ref{th:Main}.
Note that utilizing the suboptimal withincomponent density leads to a significant reduction of the asymptotic constant, as compared to
the uniform withincomponent knot distribution.

\section{Proofs}

\noindent{\it Proof of Theorem \ref{th:Main}.}
First we investigate the asy\-mpto\-tic be\-havi\-or of the approximation error
$e_{N}(\mathbf{t}) : = \norm{X(\mathbf{t})- X_N(\mathbf{t})}$ for any  $\mathbf{t}\in \dcal_{\mathbf{i}}$, $\mb{i}\in\mb{I}$, where
$\mb{I}:=\{\mathbf{i}=(i_1,\ldots,i_d)$, $0\leq i_k\leq n_k^*-1$, $k=1,\ldots,d\}$,
when the number of knots $N$ tends to infinity.
Further, we find the asymptotic form  of the IMSE
$$
e_N:=\left(\int_\dcal e_N(\mathbf{t})^2 d\mathbf{t}\right)^{1/2}
$$
for any positive continuous densities $h_1(\cdot),\ldots,h_k(\cdot)$.
 We start by observing that
\begin{eqnarray}
 e_N(\mathbf{t})^2 \!\! &=&\!\!\!\E(X(\mathbf{t})-X_N(\mathbf{t}))^2=\E(\E_{\boldsymbol\eta}(X(\mathbf{t}_{\mathbf{i}}+
               \mathbf{r}_{\mathbf{i}}\!\!\ast\!\! {\boldsymbol\eta})-X(\mathbf{t})))^2\nonumber \\
            &=&\!\!\!\E_{\boldsymbol\eta,\boldsymbol\xi} \E\left(\left(X(\mathbf{t}_{\mathbf{i}}+\mathbf{r}_{\mathbf{i}}\!\!\ast\!\!{\boldsymbol\eta})-
               X(\mathbf{t})\right)\left(X(\mathbf{t}_{\mathbf{i}}+\mathbf{r}_{\mathbf{i}}\!\!\ast\!\!{\boldsymbol\xi})-X(\mathbf{t})\right)\right)\nonumber\\
&=&\!\!\!\frac{1}{2}\E_{\boldsymbol \eta,\boldsymbol\xi}\E\left((X(\mathbf{t}_{\mathbf{i}}+\mathbf{r}_{\mathbf{i}}\!\ast\! \boldsymbol\eta)\!\!-\!\!
                      X(\mathbf{t}))^2\!\! +\!\! (X(\mathbf{t}_{\mathbf{i}}+\mathbf{r}_{\mathbf{i}}\!\ast\!\boldsymbol\xi)\!\!-\!\!X(\mathbf{t}))^2
                \!-\!(X(\mathbf{t}_{\mathbf{i}}+\mathbf{r}_{\mathbf{i}}\!\ast\!\boldsymbol\eta)\!\!-\!\!X(\mathbf{t}_{\mathbf{i}}+\mathbf{r}_{\mathbf{i}}\!\ast\!\boldsymbol \xi))^2 \right),
\label{eq:ThMain1}
\end{eqnarray}
where $\boldsymbol\xi$ is an independent copy of $\boldsymbol\eta$.
Further, the property \eqref{def:locstat} together with the uniform continuity and positiveness of local stationarity functions $c_1(\cdot),\ldots,c_k(\cdot)$ imply that
\begin{equation}\label{eq:ThMain2}
e_N(\mathbf{t})^2 =\frac{1}{2}\left(\sum_{j=1}^{k}c_j(\mathbf{t}_{\mathbf{i}})\E_{\boldsymbol\eta,\boldsymbol\xi}\left(
\norm{\mathbf{r}^j_{\mathbf{i}}\!\ast\!(\boldsymbol\eta^j-\mathbf{s}^j)}^{\alpha_j}
\!\!\!\!+\norm{\mathbf{r}^j_{\mathbf{i}}\!\ast\!(\boldsymbol\xi^j-\mathbf{s}^j)}^{\alpha_j}
\!\!\!\!-\norm{\mathbf{r}^j_{\mathbf{i}}\!\ast\!(\boldsymbol\eta^j-\boldsymbol\xi^j)}^{\alpha_j}
\right)\right)
(1+q_{N,\mb{i}}(\mathbf{t})),
\end{equation}
where $\varepsilon_N :=\max\{|q_{N,\mb{i}}(\mathbf{t})|,\mb{t}\in\dcal_{\mb{i}},\mb{i}\in\mb{I}\}=\mathrm{o}(1)$
as $N\to\infty$ \citep[cf.][]{Seleznjev2000}.
It follows from the definition and the mean (integral) value theorem that
$$
\mathbf{r}_{\mathbf{i}}=\left(\frac{1}{h_{1}^*(w_{1,i_1})n_{1}^*},\frac{1}{h_{2}^*(w_{2,i_2})n_{2}^*},\ldots,\frac{1}{h_{d}^*(w_{d,i_d})n_{d}^*}\right),
      \quad w_{j,i_j}\in [t_{j,i_j},t_{j,i_j+1}],\,j=1,\ldots,d.
$$
Denote by $\mathbf{w}_{\mathbf{i}}:=(w_{1,i_1},\ldots,w_{d,i_d})$. Now the definition of $cRS(h,\pi,\mb{l})$ implies
$$
\mathbf{r}_{\mathbf{i}}^j=\left(\frac{1}{n_j h_j(w_{L_{j-1}+1,i_{L_{j-1}+1}})},\ldots,\frac{1}{n_j h_j(w_{L_j,i_{L_j}})}\right)=\frac{1}{n_j} H_j(\mb{w}_{\mb{i}}^j),\quad j=1,\ldots,k,
$$
where $H_j(\mb{t}^j):=(1/h_j(t_{L_{j-1}+1}),\ldots,1/h_j(t_{L_{j}}))$, $j=1,\ldots,k$. Consequently,
\begin{eqnarray}
e_N(\mathbf{t})^2&=& \frac{1}{2}\Bigg(\sum_{j=1}^{k}n_j^{-\alpha_j}c_j(\mathbf{t}_{\mathbf{i}})
\E_{\eta,\xi}\Big(||H_j(\mb{w}^j_{\mb{i}})\ast(\boldsymbol\eta^{j}-\mathbf{s}^j)||^{\alpha_j}
+||H_j(\mb{w}^j_{\mb{i}})\ast(\boldsymbol\xi^{j}-\mathbf{s}^j)||^{\alpha_j}\nonumber\\
&&\qquad\qquad-||H_j(\mb{w}^j_{\mb{i}})\ast(\boldsymbol\eta^{j}-{\boldsymbol \xi}^j)||^{\alpha_j}\Big)
\Bigg)(1+\mathrm{o}(1))\quad \ninf.\nonumber
\end{eqnarray}
Applying the uniform continuity of $h(\cdot)$ yields
$$
\begin{aligned}
e_N(t)^2&=\frac{1}{2}\Bigg(\sum_{j=1}^{k}n_j^{-\alpha_j}c_j(\mb{t_i})
\E_{\eta,\xi}\Big(||H_j(\mb{t}^j_{\mb{i}})\ast(\boldsymbol\eta^{j}-\mathbf{s}^j)||^{\alpha_j}
+||H_j(\mb{t}^j_{\mb{i}})\ast(\boldsymbol\xi^{j}-\mathbf{s}^j)||^{\alpha_j}\\
&\qquad\qquad-||H_j(\mb{t}^j_{\mb{i}})\ast(\boldsymbol\eta^{j}-{\boldsymbol \xi}^j)||^{\alpha_j}\Big)\Bigg)(1+\mathrm{o}(1))\\
&=\Bigg(\sum_{j=1}^{k}n_j^{-\alpha_j}c_j(\mb{t_i})C_{\alpha_j,l_j}(\mb{s}^j;H_j(\mb{t}^j_{\mb{i}}))
\Bigg)(1+\mathrm{o}(1))\quad \ninf,
\end{aligned}
$$
where
$$
\begin{aligned}
C_{\alpha_j,l_j}(\mb{s}^j;H_j(\mb{t}^j_{\mb{i}})) :=&
\frac{1}{2}\E_{\boldsymbol\eta,\boldsymbol\xi}\Big(||H_j(\mb{t}^j_{\mb{i}})\ast(\boldsymbol\eta^{j}-\mathbf{s}^j)||^{\alpha_j}
+||H_j(\mb{t}^j_{\mb{i}})\ast(\boldsymbol\xi^{j}-\mathbf{s}^j)||^{\alpha_j}
-||H_j(\mb{t}^j_{\mb{i}})\ast(\boldsymbol\eta^{j}-{\boldsymbol \xi}^j)||^{\alpha_j}\Big)\\
=&\norm{B_{\alpha_j,l_j}(H_j(\mb{t}^j_{\mb{i}}) \ast \mb{s}^j)-\E_{\boldsymbol\eta} B_{\alpha_j,l_j}(H_j(\mb{t}^j_{\mb{i}})\!\ast\! \boldsymbol\eta^j)}^2_2.
\end{aligned}
$$
\\%
Let $\dcal_{\mb{i}}=\dcal_{\mb{i}}^1\times\cdots\times\dcal_{\mb{i}}^k$ and denote by $|\dcal_{\mb{i}}|$
the volume of hyperrectangle $\dcal_{\mb{i}}$. Then
$$
\begin{aligned}
e_N^2   &=\sum_{\mb{i}\in\mb{I}} \int_{\dcal_{\mb{i}}}e_N(\mb{t})^2 d\mb{t} =
    \Bigg(\sum_{\mb{i}\in\mb{I}} \int_{\dcal_{\mb{i}}} \sum_{j=1}^{k}n_j^{-\alpha_j}c_j(\mb{t_i})C_{\alpha_j,l_j}(\mb{s}^j;H_j(\mb{t}^j_{\mb{i}})) d\mb{t}\Bigg)(1+\mathrm{o}(1))\\
    &=\Bigg(\sum_{\mb{i}\in\mb{I}} \sum_{j=1}^{k}n_j^{-\alpha_j}c_j(\mb{t_i}) \int_{\dcal^j}C_{\alpha_j,l_j}(\mb{s}^j;H_j(\mb{t}^j_{\mb{i}}))d\mb{s}^j|\dcal_{\mb{i}}|\Bigg)(1+\mathrm{o}(1))\\
    &=\Bigg(\sum_{j=1}^{k}n_j^{-\alpha_j} \sum_{\mb{i}\in\mb{I}} c_j(\mb{t_i}) b_{\alpha_j,l_j}(H_j(\mb{t}^j_{\mb{i}}))|\dcal_{\mb{i}}|\Bigg)(1+\mathrm{o}(1)) \ninf. \\
\end{aligned}
$$
Now the Riemann integrability of the functions $c_j(\mb{t}) b_{\alpha_j,l_j}(H_j(\mb{t}^j))$, $j=1,\ldots,k$, gives
$$
\begin{aligned}
e_N^2&=\Bigg(\sum_{j=1}^{k}n_j^{-\alpha_j}\int_\dcal c_j(\mb{t}) b_{\alpha_j,l_j}(H_j(\mb{t}^j))d\mb{t}\Bigg)(1+\mathrm{o}(1))=\Bigg(\sum_{j=1}^{k}\frac{v_j}{n_j^{\alpha_j}}\Bigg)(1+\mathrm{o}(1))\quad \ninf.
\end{aligned}
$$
Note that for any $\mb{u}\in \mathbb{R}^{m}_{+}$,  $b_{\beta,m}(\mb{u})>0$, otherwise the fractional Brownian field is degenerated \citep[cf.][]{Seleznjev2000}. Consequently,
$v_j>0$, $j=1,\ldots,k$. This completes the proof.
\bigskip

\noindent{\it Proof of Theorem \ref{Th:DimOpt}.} Note that by the inequality for the arithmetic and geometric means,
$$
\frac{1}{k}\sum_{j=1}^{k}\frac{v_j}{n_j^{\alpha_j}}\geq\left(\prod_{j=1}^{k}\frac{v_j}{n_j^{\alpha_j}}\right)^{1/k}
$$
with equality if only if
$$
\nu^{-1}=\frac{v_j}{n_j^{\alpha_j}},\quad j=1,\ldots,k.
$$
Hence, the equality is attained for $\tilde n_j=(\nu v_j)^{1/\alpha_j}$, $j=1,\ldots,k$. Let
\begin{equation}\label{eq:ThDimOpt1}
n_j = \left\lceil\tilde n_j\right\rceil\sim\left({\nu}{v_j}\right)^{1/\alpha_j} \mbox{ as } N\to\infty.
\end{equation}
The total number of observations satisfies
$$
N=(n_1^*+1)\cdots(n_d^*+1)\sim\prod_{i=1}^{d}n_i^*=\prod_{j=1}^{k}n_j^{l_j}=M\quad \mbox{ as } N\to\infty.
$$
This implies that for the asymptotically optimal intercomponent knot distribution
$$
N\sim M \sim \nu^{1/\rho}\prod_{j=1}^{k}{v_j^{l_j/\alpha_j}},
$$
and therefore,
$$
\nu\sim {N^{\rho}}{\kappa^{-\rho}} \mbox{ as } N\to\infty.
$$
By equation \eqref{eq:ThDimOpt1}, the asymptotically optimal intercomponent knot distribution is
$$
n_j\sim\frac{N^{\rho/\alpha_j} v_j^{1/\alpha_j} }{\kappa^{\rho/\alpha_j}}\mbox{ as } N\to\infty,\quad j=1,\ldots,k.
$$
Moreover, with such chosen knot distribution, the equality in \eqref{eq:Th2} is attained asymptotically. This completes the proof.
\\\bigskip

\noindent{\it Proof of Proposition \ref{Prop:OneDim}.} The proof is a straightforward implication of the assumptions and
equation \eqref{eq:ThMain2}. The exact constant and the expression for the optimal density are
due to \citet{Seleznjev2000}.
\bigskip

\noindent{\it Proof of Proposition \ref{th:Suboptimal}.} The first steps of the proof repeat those of Theorem \ref{th:Main}. By \eqref{eq:ThMain2}, we have
$$
\begin{aligned}
e_N(\mathbf{t})^2 &=\frac{1}{2}\Bigg(\sum_{j=1}^{k}c_j(\mathbf{t}_{\mathbf{i}})\E_{\boldsymbol\eta,\boldsymbol\xi}\left(
\norm{\mathbf{r}^j_{\mathbf{i}}\!\ast\!(\boldsymbol\eta^j-\mathbf{s}^j)}^{\alpha_j}
+\norm{\mathbf{r}^j_{\mathbf{i}}\!\ast\!(\boldsymbol\xi^j-\mathbf{s}^j)}^{\alpha_j}
-\norm{\mathbf{r}^j_{\mathbf{i}}\!\ast\!(\boldsymbol\eta^j-\boldsymbol\xi^j)}^{\alpha_j}
\right)\Bigg)(1+\mathrm{o}(1))\\
&\leq\frac{1}{2}\Bigg( \sum_{j=1}^{k}c_j(\mathbf{t}_{\mathbf{i}})\E_{\boldsymbol\eta,\boldsymbol\xi}\left(
\norm{\mathbf{r}^j_{\mathbf{i}}\!\ast\!(\boldsymbol\eta^j-\mathbf{s}^j)}^{\alpha_j}
+\norm{\mathbf{r}^j_{\mathbf{i}}\!\ast\!(\boldsymbol\xi^j-\mathbf{s}^j)}^{\alpha_j}
\right)\Bigg)(1+\mathrm{o}(1))\\
&=\Bigg(\sum_{j=1}^{k}c_j(\mathbf{t}_{\mathbf{i}})\E_{\boldsymbol\eta}\left(
\norm{\mathbf{r}^j_{\mathbf{i}}\!\ast\!(\boldsymbol\eta^j-\mathbf{s}^j)}^{\alpha_j}\right)\Bigg)(1+\mathrm{o}(1)) \ninf.
\end{aligned}
$$
For any nonnegative numbers $a_1,\ldots,a_k$ and any $\alpha\in \mathbb{R_+}$, the inequality
\begin{equation}\label{PowerIneq}
\left(\sum_{i=1}^{k} a_i\right)^\alpha \leq k^\alpha\sum_{i=1}^{k}a_i^\alpha
\end{equation}
 holds, and consequently,
$$
\begin{aligned}
e_N(\mb{t})^2 &\leq\Bigg( \sum_{j=1}^{k}c_j(\mb{t_i}) l_j^{\alpha_j/2}\sum_{m=L_{j-1}+1}^{L_j}\E_{\boldsymbol\eta}(r_{\mb{i},m}|\eta_{m}-s_m|)^{\alpha_j}\Bigg)(1+\mathrm{o}(1))\\
         &=   \Bigg(   \sum_{j=1}^{k}c_j(\mb{t_i}) l_j^{\alpha_j/2}\sum_{m=L_{j-1}+1}^{L_j}r_{\mb{i},m}^{\alpha_j}\left((1-s_m)^{\alpha_j}s_m+(1-s_m)s_m^{\alpha_j}\right)\Bigg) (1+\mathrm{o}(1)).\\
\end{aligned}
$$
By the mean value theorem and the uniform continuity of withincomponent densities, we obtain
$$
\begin{aligned}
e_N(\mb{t})^2 &\leq\Bigg( \sum_{j=1}^{k}c_j(\mb{t_i}) l_j^{\alpha_j/2}n_j^{-\alpha_j}\!\!\!\sum_{m=L_{j-1}+1}^{L_j}(h_j(t_{\mb{i},m}))^{-\alpha_j}\left((1-s_m)^{\alpha_j}s_m+(1-s_m)s_m^{\alpha_j}\right)\Bigg)(1+\mathrm{o}(1))\ninf.\\
\end{aligned}
$$
Proceeding now to the calculation of the IMSE, we get
$$
\begin{aligned}
e_N^2&=\sum_{\mb{i}\in\mb{I}} \int_{\dcal_{\mb{i}}}e_N(\mb{t})^2 dt
\leq \Bigg(\sum_{\mb{i}\in\mb{I}}  \sum_{j=1}^{k}c_j(\mb{t_i}) l_j^{\alpha_j/2}n_j^{-\alpha_j}  \sum_{m=L_{j-1}+1}^{L_j}(h_j(t_{\mb{i},m}))^{-\alpha_j}
\frac{2}{(\alpha_j+1)(\alpha_j+2)}|\dcal_{\mb{i}}|\Bigg)(1+\mathrm{o}(1)),\\
\end{aligned}
$$
where
$$
\frac{2}{(\alpha_j+1)(\alpha_j+2)}=\int_0^1\left((1-s)^{\alpha_j}s+(1-s)s^{\alpha_j}\right)ds.
$$
Now the Riemann integrability of
$c_j(\mb{t}) h_j(t_m)^{-\alpha_j}$, $j=1,\ldots,k$, together with the definition of integrated local stationarity functions imply that
$$
\begin{aligned}
e_N^2&\leq\Bigg(\sum_{j=1}^{k}\frac{1}{n_j^{\alpha_j}}{l_j^{\alpha_j/2}}\left(a_{\alpha_j}+\frac{1}{6}\right)\sum_{m=L_{j-1}+1}^{L_j}\sum_{\mb{i}\in\mb{I}} c_j(\mb{t_i}) (h_j(t_{\mb{i},m}))^{-\alpha_j} |\dcal_{\mb{i}}|\Bigg)(1+\mathrm{o}(1))\\
&=\Bigg(\sum_{j=1}^{k}\frac{1}{n_j^{\alpha_j}}{l_j^{\alpha_j/2}}\left(a_{\alpha_j}+\frac{1}{6}\right)\sum_{m=L_{j-1}+1}^{L_j}\int_\dcal c_j(\mb{t}) h_j(t_m)^{-\alpha_j} dt\Bigg)(1+\mathrm{o}(1))\\
&=\Bigg(\sum_{j=1}^{k}\frac{1}{n_j^{\alpha_j}}{l_j^{1+\alpha_j/2}}\left(a_{\alpha_j}+\frac{1}{6}\right)\int_0^1 C_j(t_{L_j}) h_j(t_{L_j})^{-\alpha_j} dt_{L_j}\Bigg)(1+\mathrm{o}(1)) \ninf.
\end{aligned}
$$
The expression for the suboptimal density is due to \citet{Seleznjev2000}. This completes the proof.
\bigskip

\noindent{\it Proof of Proposition \ref{th:HoldType}.} We start by proving $(i)$.
Let $X\in\ccal^{\boldsymbol\alpha}_\mb{l}([0,1]^d,C)$ and consider $\mb{t}\in\dcal_\mb{i}$, $\mb{i}\in\mb{I}$. Applying the H\"{o}lder condition \eqref{def:hcont} to equation \eqref{eq:ThMain1} yields
$$
\begin{aligned}
  e_N(\mathbf{t})^2     &=\frac{1}{2}\E_{\boldsymbol \eta,\boldsymbol\xi}\E\left((X(\mathbf{t}_{\mathbf{i}}+\mathbf{r}_{\mathbf{i}}\!\ast\! \boldsymbol\eta)-
                      X(\mathbf{t}))^2\!\! +\!\! (X(\mathbf{t}_{\mathbf{i}}+\mathbf{r}_{\mathbf{i}}\!\ast\!\boldsymbol\xi)-X(\mathbf{t}))^2
                \!-\!(X(\mathbf{t}_{\mathbf{i}}+\mathbf{r}_{\mathbf{i}}\!\ast\!\boldsymbol\eta)-X(\mathbf{t}_{\mathbf{i}}+\mathbf{r}_{\mathbf{i}}\!\ast\!\boldsymbol \xi))^2 \right) \\
&\leq C \E_{\boldsymbol\eta}\norm{\mb{r_i}\ast\boldsymbol\eta}_{\boldsymbol\alpha}=C \E_{\boldsymbol\eta}\sum_{j=1}^k\norm{\mb{r_i}^j\ast\boldsymbol\eta^j}^{\alpha_j}\leq
C\sum_{j=1}^{k}l_j^{\alpha_j/2}\!\!\!\!\sum_{m=L_{j-1}+1}^{L_j}\E_{\boldsymbol\eta}\left({r}_{\mb{i},m}|\eta_m-s_m|\right)^{\alpha_j},\\
\end{aligned}
$$
where the last inequality follows from \eqref{PowerIneq}.
Furthermore, since $\max_{s\in[0,1]}\left((1-s)^{\alpha_j} s+(1-s) s^{\alpha_j}\right)=2^{-\alpha_j}$, we obtain
$$
\begin{aligned}
  e_N(\mathbf{t})^2 &\leq
C\sum_{j=1}^{k} l_j^{\alpha_j/2}\!\!\!\!\sum_{m=L_{j-1}+1}^{L_j}r_{\mb{i},m}^{\alpha_j} \left((1-s_m)^{\alpha_j} s_m+(1-s_m) s_m^{\alpha_j} \right)
\leq
 \sum_{j=1}^{k} 2^{-\alpha_j} l_j^{\alpha_j/2}\!\!\!\!\sum_{m=L_{j-1}+1}^{L_j}r_{\mb{i},m}^{\alpha_j}.
\end{aligned}
$$
By the regularity of the generating densities, we have that $r_{\mb{i},m}\leq 1/(n_m^* \min_{s\in[0,1]}h^*_m(s) ),\mb{i}\in\mb{I},m=1,\ldots,d$.
Moreover, the definition of $cRS(h,\pi,\mb{l})$ implies the following uniform bound for the squared approximation accuracy
$$
\norm{X-X_N}^2_\infty=\max_{\mb{t}\in\dcal}e_N^2(\mb{t})\leq
C\sum_{j=1}^{k} 2^{-\alpha_j}l_j^{1+\alpha_j/2}\left(\frac{D_j}{n_j}\right)^{\alpha_j},\\
$$
with $D_j=1/\min_{s\in[0,1]}h_j(s)$, $j=1,\ldots,k$. Finally, we obtain the required assertion
$$
\norm{X-X_N}_\infty\leq \sqrt{C}\sum_{j=1}^{k} \frac{c_j}{n_j^{\alpha_j/2}},
$$
where $c_j^2:=2^{-\alpha_j} l_j^{1+\alpha_j/2} D_j^{\alpha_j}>0$, $j=1,\ldots,k$.

For the smooth case, we use the multivariate Taylor formula to obtain the following
representation of the deviation field
$$
\begin{aligned}
\delta_n(\mb{t}):=X(\mb{t})-X_N(\mb{t})=\E_{\boldsymbol\eta}\left(\int_0^1\sum_{j=1}^d X_j'(\mb{t_{i}}+u\, \mb{r_i}\ast(\boldsymbol\eta-\mb{s})){r}_{\mb{i},j}(\eta_j-s_j)du \right),
\quad \mb{t}\in \dcal_\mb{i},\,\mb{t}=\mb{t_i}+\mb{s}\ast \mb{r_{i}},
\end{aligned}
$$
where $\boldsymbol\eta=(\eta_1,\ldots,\eta_d)$ and $\eta_1,\ldots,\eta_d$ are independent Bernoulli random variables, $\eta_j\in Be(s_j),\,j=1,\ldots,d$.
Introducing an auxiliary uniform random variable \mbox{$U\in\mathscr{U}(0,1)$} we get
$$
\begin{aligned}
\delta_n(\mb{t})&= \sum_{j=1}^{d}\E_{\boldsymbol\eta,U}\left( X_j'(\mb{t_i}+U(\boldsymbol\eta-\mb{s})\!\ast\! \mb{r_i}){r}_{\mb{i},j} (\eta_j-s_j)\right)\\
       &= \sum_{j=1}^{d} \E_{\boldsymbol\eta,U}  \Big(  X_j'(t_{\mb{i},1}+U(\eta_1-s_1)r_{\mb{i},1},\ldots,t_{\mb{i},j}+U(\eta_j-s_j)r_{\mb{i},j},\ldots,t_{\mb{i},d}+U(\eta_d-s_d)r_{\mb{i},d})\\
       &\qquad\qquad\qquad-X_j'(t_{\mb{i},1}+U(\eta_1-s_1)r_{\mb{i},1},\ldots,t_{\mb{i},j},\ldots,t_{\mb{i},d}+U(\eta_d-s_d)r_{\mb{i},d})\Big)(\eta_j-s_j),
\end{aligned}
$$
since for any $j=1,\ldots,d$,
$$
\begin{aligned}
&\E_{\boldsymbol\eta}(X_j'(t_{\mb{i},1}+U(\eta_1-s_1)r_{\mb{i},1},\ldots,t_{\mb{i},j},\ldots,t_{\mb{i},d}+U(\eta_d-s_d)r_{\mb{i},d})(\eta_j-s_j))\\
&=\E_{\eta_1,\ldots,\eta_{j-1},\eta_{j+1},\ldots,\eta_d}
(X_j'(t_{\mb{i},1}+U(\eta_1-s_1)r_{\mb{i},1},\ldots,t_{\mb{i},j},\ldots,t_{\mb{i},d}+U(\eta_d-s_d)r_{\mb{i},d})\E_{\eta_j}(\eta_j-s_j))=0.
\end{aligned}
$$
The triangle inequality and the condition \eqref{def:hdiff} imply that
$$
e_N(\mb{t})\leq \sum_{j=1}^{d} \sqrt{C}\, V_j\, r_{\mb{i},j}^{1+\alpha_j/2},
$$
for some positive constants $V_j,j=1,\ldots,d$. Analogously to $(i)$, the required assertion follows from the regularity of the generating densities and the definition of $cRS(h,\pi,\mb{l})$.
This completes the proof.\\
\bigskip

\noindent{\bf Acknowledgments}
\smallskip\hfill

\noindent
The second author is partly supported by the Swedish Research Council grant 2009-4489 and the project "Digital Zoo" funded by the European Regional Development Fund.

\bibliographystyle{model2-names}

\bibliography{ApproxFields}
\end{document}